\documentclass[10 pt]{amsart}
\usepackage{amsfonts}
\usepackage{ifthen}
\usepackage{amsthm}
\usepackage{amsmath}
\usepackage{graphicx}
\usepackage{amscd,amssymb,amsthm}
\usepackage{graphicx}
\usepackage{epstopdf}
\usepackage{hyperref}

\newcounter{minutes}
\divide\time by 60
\newcounter{hours}
\multiply\time by 60 \addtocounter{minutes}{-\time}

\newtheorem{lemma}{Lemma}
\newtheorem{theorem}{Theorem}
\newtheorem{corollary}{Corollary}

\newcommand{\real}{\operatorname{Re}}


\begin{document}
\title{Radii of starlikeness and convexity of Wright functions}

\author[\'A. Baricz]{\'Arp\'ad Baricz$^{\bigstar}$}
\address{Department of Economics, Babe\c{s}-Bolyai University, Cluj-Napoca 400591, Romania}
\address{Institute of Applied Mathematics, \'Obuda University, 1034	Budapest, Hungary}
\email{bariczocsi@yahoo.com}

\author[E. Toklu]{Evr{\.i}m Toklu}
\address{Department of Mathematics, Faculty of Science, A\u{g}r{\i} {\.I}brah{\.I}m \c{C}e\c{c}en University, 04100 A\u{g}r{\i}, Turkey}
\email{evrimtoklu@gmail.com}

\author[E. Kad{\i}o{\u{g}}lu]{Ekrem Kad{\i}o{\u{g}}lu}
\address{Department of Mathematics, Faculty of Science, Atat\"urk University, 25240 Erzurum, Turkey}
\email{ekadioglu@atauni.edu.tr}

\thanks{$^{\bigstar}$The research of \'A. Baricz was supported by a research grant of the Babe\c{s}-Bolyai University for young researchers with project number GTC-31777.}

\def\thefootnote{}
\footnotetext{ \texttt{File:~\jobname .tex,
		printed: \number\year-0\number\month-\number\day,
		\thehours.\ifnum\theminutes<10{0}\fi\theminutes}
} \makeatletter\def\thefootnote{\@arabic\c@footnote}\makeatother

\maketitle
\begin{abstract}
In this paper our aim is to find the radii of starlikeness and convexity of the normalized Wright functions for three different kind of normalization. The key tools in the proof of our main results are the Mittag-Leffler expansion for Wright function and properties of real zeros of Wright function and its derivative. In addition, by using the Euler-Rayleigh inequalities we obtain some tight lower and upper bounds for the radii of starlikeness and convexity of order zero for  the normalized Wright functions. The main results of the paper are natural extensions of some known results on classical Bessel functions of the first kind. Some open problems are also proposed, which may be of interest for further research.

\end{abstract}

\section{\bf Introduction}
Special functions are indispensable in many branches of mathematics and applied mathematics. Geometric properties of some special functions were examined by many authors (see \cite{abo,aby,bdoy,Ba0,Ba2,Ba1,Ba3,basz2,L1,bçdt,mustafa,p1,p2,raza,szasz,Sza}). However, its origins can be traced to Brown \cite{Brown0} (see also \cite{Brown1,Brown2}), to Kreyszig and Todd \cite{KT} and to Wilf \cite{wilf}. Recently, there has been a vivid interest on geometric properties of special functions such as Bessel, Struve, Lommel functions of the first kind and regular Coulomb wave functions. The first author and his collaborators examined in details the determination of the radii of starlikeness and convexity of some normalized forms of these special functions, see \cite{abo,aby,bdoy,Ba0,Ba2,Ba1,Ba3,basz2,L1,bçdt} and the references therein. Moreover, one of the most important things which we have learned in these studies is that the radii of univalence, starlikeness and convexity are obtained as solutions of some transcendental equations and the obtained radii satisfy some interesting inequalities. In addition, in view of these studies, we know that the radii of univalence of some normalized Bessel, Struve, Lommel and regular Coulomb wave functions coincide with the radii of starlikeness of the these functions. The positive zeros of Bessel, Struve, Lommel functions of the first kind and regular Coulomb wave functions and the Laguerre-P\'olya class of real entire functions played an important role in these papers. Motivated by the above series of papers on geometric properties of special functions, in this paper our aim is to present some similar results for the normalized forms of the Wright function which has important applications in different areas of mathematics. In this paper, we are mainly focused on the determination of the radii of starlikeness and convexity of the normalized Wright functions. Furthermore, our aim is also to give some lower and upper bounds for the radii of starlikeness and convexity of order zero by using some Euler-Rayleigh inequalities for the smallest positive zero of some transcendental equations (for more details on such kind of inequalities we refer to \cite{ismail}). The paper is arranged as follows: the rest of this section is devoted to some basic definitions, which are needed for the proof of our main results. Section 2 is separated into four subsections: the first subsection is dedicated to the radii of starlikeness of normalized Wright functions. Also, at the end of this subsection there are given lower and upper bounds for radii of starlikeness of order zero. The second subsection contains the study of the radii of convexity of normalized Wright functions, and in its last part lower and upper bounds for radii of convexity of order zero for some normalized Wright functions are given. The third subsection contains some particular cases of the main results in terms of the classical Bessel functions of the first kind. In the fourth subsection some open problems are stated, which may be of interest for further research. The proofs of the main results are presented in the third section.

Before to start the presentation of the results we would like to state some basic definitions. For $r>0$ we denote by $\mathbb{D}_r=\left\{z\in\mathbb{C}: |z|<r\right\}$ the open disk of radius $r$
centered at the origin. Let $f:\mathbb{D}_r\to\mathbb{C}$ be the function defined by
\begin{equation}
f(z)=z+\sum_{n\geq 2}a_{n}z^{n},  \label{eq0}
\end{equation}
where $r$ is less or equal than the radius of convergence of the above power series. Let $\mathcal{A}$ be the class of analytic functions of the form \eqref{eq0}, that is, normalized by the conditions $f(0)=f^{\prime}(0)-1=0.$ The
function $f,$ defined by \eqref{eq0}, is called starlike in $\mathbb{D}_r$
if $f$ is univalent in $\mathbb{D}_r$, and the image domain $f(\mathbb{D}_r)$ is a
starlike domain in $\mathbb{C}$ with respect to the origin (see \cite{Dur} for more details). Analytically, the function $f$ is starlike in
$\mathbb{D}_r$ if and only if $$\real\left( \frac{zf^{\prime }(z)}{f(z)}
\right) >0 \quad \mbox{for all}\ \ z\in
\mathbb{D}_r.$$ For $\alpha \in [0,1)$ we say that
the function $f$ is starlike of order $\alpha $ in $\mathbb{D}_r$ if and
only if $$\real\left( \frac{zf^{\prime }(z)}{f(z)}\right) >\alpha \quad \mbox{for all}\ \ z\in
\mathbb{D}_r.$$
The real number
\begin{equation*}
r_{\alpha }^{\star}(f)=\sup \left\{ r>0 \left|\real
\left( \frac{zf^{\prime }(z)}{f(z)}\right)  >\alpha \;\text{for all }z\in
\mathbb{D}_r\right.\right\}
\end{equation*}
is called the radius of starlikeness of order $\alpha$ of the function $f$. Note that $
r^{\star}(f)=r_{0}^{\star}(f)$ is in fact the largest radius such that the image
region $f(\mathbb{D}_{r^{\star}(f)})$ is a starlike domain with respect to
the origin.

The function $f,$ defined by \eqref{eq0}, is convex in the disk $\mathbb{D}_r$ if $f$
is univalent in $\mathbb{D}_r$, and the image domain $f(\mathbb{D}_r)$ is a
convex domain in $\mathbb{C}.$ Analytically, the function $f$ is convex in $\mathbb{D}_r$ if and only
if $$\real\left(  1+\frac{zf^{\prime \prime }(z)}{f^{\prime }(z)}\right)
>0  \quad \mbox{for all}\ \ z\in
\mathbb{D}_r.$$
For $\alpha \in[0,1)$ we say that the
function $f$ is convex of order $\alpha $ in $\mathbb{D}_r$ if and only if $$
\real\left( 1+\frac{zf^{\prime \prime }(z)}{f^{\prime }(z)}\right)
>\alpha \quad \mbox{for all}\ \ z\in
\mathbb{D}_r.$$ The radius of convexity of order $\alpha $
of the function $f$ is defined by the real number%
\begin{equation*}
r_{\alpha }^{c}(f)=\sup \left\{ r>0 \left|\real\left( 1+
\frac{zf^{\prime \prime }(z)}{f^{\prime }(z)}\right) >\alpha \;\text{for all }
z\in\mathbb{D}_r\right.\right\} .
\end{equation*}
Note that $r^{c}(f)=r_{0}^{c}(f)$ is the largest radius such that the image region $f(\mathbb{D}_{r^{c}(f)})$ is a convex domain. 

Finally, we recall that a real entire function $q$ belongs to the  Laguerre-P\'{o}lya class $\mathcal{LP}$ if it can be represented in the form $$q(x)=cx^{m}e^{-ax^2+bx}\prod_{n\geq1}\left(1+\frac{x}{x_n}\right)e^{-\frac{x}{x_n}},$$ with $c,b,x_n\in\mathbb{R}, a\geq0, m\in\mathbb{N}_0$ and $\sum1/{x_n}^2<\infty.$ We note that the class $\mathcal{LP}$ consists of entire functions which are uniform limits on the compact sets of the complex plane of polynomials with only real zeros. For more details on the class $\mathcal{LP}$ we refer to \cite[p. 703]{DC} and to the references therein.

\section{\bf The radii of starlikeness and convexity of normalized Wright functions}

In this section, we will investigate the generalized Bessel function
$$\phi(\rho,\beta,z)=\sum_{n\geq0}\frac{z^n}{n!\Gamma(n\rho+\beta)},$$
where $\rho>-1$ and $z,\beta\in\mathbb{C},$ named after E.M. Wright. This function was introduced by Wright for $\rho>0$ in connection with his investigations on the asymptotic theory of partitions \cite{wright}, see also \cite{mainardi} for further details. Furthermore, it is important to mention that the Wright function is an entire function of $z$ for $\rho>-1$, consequently, as we will see in some parts of our paper, some properties of the general theory of entire functions can be applied.

The following Lemma, which we believe it is of independent interest, plays an important role in the proof of our main results.

\begin{lemma}\label{lemma}
If $\rho>0$ and $\beta>0,$ then the function $z\mapsto\lambda_{\rho,\beta}(z)=\phi(\rho,\beta,-z^2)$ has infinitely many zeros which are all real. Denoting by $\lambda_{\rho,\beta,n}$ the $n$th positive zero of $\phi(\rho,\beta,-z^2),$ under the same conditions the Weierstrassian decomposition 
$$\Gamma(\beta)\phi(\rho,\beta,-z^2)=\prod_{n\geq1}\left(1-\frac{z^2}{\lambda^2_{\rho,\beta,n}}\right)$$
is valid, and this product is uniformly convergent on compact subsets of the complex plane. Moreover, if we denote by $\zeta_{\rho,\beta,n}'$ the $n$th positive zero of $\Psi_{\rho,\beta}',$ where $\Psi_{\rho,\beta}(z)=z^{\beta}\lambda_{\rho,\beta}(z),$ then the positive zeros of $\lambda_{\rho,\beta}$ (or $\Psi_{\rho,\beta}$) are interlaced with those of $\Psi_{\rho,\beta}'$. In other words, the zeros satisfy the chain of inequalities
$$\zeta^{\prime}_{\rho,\beta,1}<\lambda_{\rho,\beta,1}<\zeta^{\prime}_{\rho,\beta,2}<\zeta_{\rho,\beta,2}<{\dots}.$$
\end{lemma}
Observe that the function $z\mapsto\phi(\rho,\beta,-z^2)$ does not belong to $\mathcal{A}$, and thus first we perform some natural normalization. We define three functions originating from $\phi(\rho,\beta,\cdot)$:
$$f_{\rho,\beta}(z)=\left(z^{\beta}\Gamma(\beta)\phi(\rho,\beta,-z^2)\right)^{\frac{1}{\beta}},$$
$$g_{\rho,\beta}(z)=z\Gamma(\beta)\phi(\rho,\beta,-z^2),$$
$$h_{\rho,\beta}(z)=z\Gamma(\beta)\phi(\rho,\beta,-z).$$
Obviously these functions belong to the class $\mathcal{A}$. Of course, there exist infinitely many other normalization, the main motivation
to consider the above ones is the fact that their particular cases in terms of Bessel functions appear in literature or are similar to the studied normalization in the literature.

\subsection{\bf The radii of starlikeness of order $\alpha$ of the functions $f_{\rho,\beta}$, $g_{\rho,\beta}$ and $h_{\rho,\beta}$}

 In this subsection our aim is to present some results for the radii of starlikeness of normalized Wright functions  $f_{\rho,\beta}$, $g_{\rho,\beta}$ and $h_{\rho,\beta}$. We will see in this subsection that the radii of starlikeness of order $\alpha$ of the normalized Wright functions are actually solutions of some transcendental equations. Moreover, we will also find lower and upper bounds for the radii of starlikeness of order zero.

Our first main result is the following theorem.

\begin{theorem}\label{th1}
Let $\rho>0$ and $\beta>0$.
\begin{enumerate}
\item[\bf a.] The radius of starlikeness of $f_{\rho,\beta}$ is $r^{\star}_{\alpha}\left(f_{\rho,\beta}\right)=x_{\rho,\beta,1} $, where $x_{\rho,\beta,1}$ is the smallest positive zero of the transcendental equation $$r\lambda^{\prime}_{\rho,\beta}(r)-(\alpha-1)\beta\lambda_{\rho,\beta}(r)=0.$$
\item[\bf b.] The radius of starlikeness of $g_{\rho,\beta}$ is $r^{\star}_{\alpha}\left( g_{\rho,\beta}\right)=y_{\rho,\beta,1} $, where $y_{\rho,\beta,1}$ is the smallest positive zero of the transcendental equation $$r\lambda^{\prime}_{\rho,\beta}(r)-(\alpha-1)\lambda_{\rho,\beta}(r)=0.$$
\item[\bf c.] The radius of starlikeness of $h_{\rho,\beta}$ is $r^{\star}_{\alpha}\left( h_{\rho,\beta}\right)=z_{\rho,\beta,1} $, where $z_{\rho,\beta,1}$ is the smallest positive zero of the transcendental equation
$$\sqrt{r}\lambda^{\prime}_{\rho,\beta}(\sqrt{r})-2(\alpha-1)\lambda_{\rho,\beta}(\sqrt{r})=0.$$
\end{enumerate}
\end{theorem}

The following theorems provide some tight lower and upper bounds for the radii of starlikeness of the functions considered in the above theorems. In these theorems for simplicity we use the notation $$\Delta_{a,b}(\rho,\beta)=a\Gamma(\beta)\Gamma(2\rho+\beta)-b\Gamma^2(\rho+\beta),$$
and we mention that the positivity of this expression for $a>b>0$ and $\rho,\beta>0$ is guaranteed by the log-convexity of the Euler gamma function.

\begin{theorem}\label{th2}
For $\rho,\beta>0$ the radius of starlikeness $r^{\star}(f_{\rho,\beta})$ satisfies
$$\sqrt{\frac{\Gamma(\rho+\beta)}{(\beta+2)\Gamma(\beta)}}<r^{\star}(f_{\rho,\beta})<
\sqrt{\frac{\beta(\beta+2)\Gamma(\rho+\beta)\Gamma(2\rho+\beta)}{\Delta_{(\beta+2)^{2},\beta+4}(\rho,\beta)}},$$
$$\sqrt[4]{\frac{\beta\Gamma^{2}(\rho+\beta)\Gamma(2\rho+\beta)}{\Gamma(\beta)\Delta_{(\beta+2)^{2},\beta+4}(\rho,\beta) }}<r^{\star}(f_{\rho,\beta})<\sqrt{{\displaystyle{2\beta\Gamma(\rho+\beta)\Gamma(3\rho+\beta)
\Delta_{(\beta+2)^{2},\beta+4}(\rho,\beta)}\over\displaystyle{\beta(\beta+6)\Gamma^{3}(\rho+\beta)\Gamma(2\rho+\beta)+
\Xi_{\rho,\beta}}}},$$
where
$$\Xi_{\rho,\beta}=(\beta+2)^2\Gamma(\beta)\Gamma(3\rho+\beta)\Delta_{2(\beta+2),\beta+4}(\rho,\beta).$$
\end{theorem}

\begin{theorem}\label{th3}
For $\rho,\beta>0$ the radius of starlikeness $r^{\star}(g_{\rho,\beta})$ satisfies
$$\sqrt{\frac{\Gamma(\rho+\beta)}{3\Gamma(\beta)}}<r^\star(g_{\rho,\beta})<
\sqrt{\frac{3\Gamma(\rho+\beta)\Gamma(2\rho+\beta)}{\Delta_{9,5}(\rho,\beta)}},$$
$$\sqrt[4]{\frac{\Gamma^{2}(\rho+\beta)\Gamma(2\rho+\beta)}{\Gamma(\beta)\Delta_{9,5}(\rho,\beta)}}<r^\star(g_{\rho,\beta})<
\sqrt{\frac{2\Gamma(\rho+\beta)\Gamma(3\rho+\beta)\Delta_{9,5}(\rho,\beta)}{9\Gamma(\beta)\Gamma(3\rho+\beta)\Delta_{6,5}(\rho,\beta) +7\Gamma^{3}(\rho+\beta)\Gamma(2\rho+\beta)}}.$$
\end{theorem}

\begin{theorem}\label{th4}
For $\rho,\beta>0$ the radius of starlikeness $r^{\star}(h_{\rho,\beta})$ satisfies
$$\frac{\Gamma(\rho+\beta)}{2\Gamma(\beta)}<r^\star(h_{\rho,\beta})<
\frac{2\Gamma(\rho+\beta)\Gamma(2\rho+\beta)}{\Delta_{4,3}(\rho,\beta)},$$
$$\sqrt{\frac{\Gamma^{2}(\rho+\beta)\Gamma(2\rho+\beta)}{\Gamma(\beta)\Delta_{4,3}(\rho,\beta)}}
<r^\star(h_{\rho,\beta})<\frac{\Gamma(\rho+\beta)\Gamma(3\rho+\beta)
\Delta_{4,3}(\rho,\beta)}{\Gamma(\beta)\Gamma(3\rho+\beta)\Delta_{8,9}(\rho,\beta)+2\Gamma^{3}(\rho+\beta)\Gamma(2\rho+\beta)}.$$
\end{theorem}

\subsection{The radii of convexity of order $\alpha$ of the functions $f_{\rho,\beta}$, $g_{\rho,\beta}$ and $h_{\rho,\beta}$ }
In this subsection we present the radii of convexity of order $\alpha$ for the functions  $f_{\rho,\beta}$, $g_{\rho,\beta}$ and $h_{\rho,\beta}$. In addition, we find tight lower and upper bounds of order zero for the functions $g_{\rho,\beta}$ and $h_{\rho,\beta}$.

\begin{theorem}\label{th5}
Let $\rho>0,$ $\beta>0$ and $\alpha\in[0,1).$
\begin{enumerate}
\item[\bf a.] The radius of convexity of order $\alpha$ of $f_{\rho,\beta}$ is the smallest positive root of
$$1+\frac{r\Psi^{\prime \prime}_{\rho,\beta}(r)}{\Psi^{\prime}_{\rho,\beta}(r)}+\left( \frac{1}{\beta}-1\right) \frac{r\Psi^{\prime}_{\rho,\beta}(r)}{\Psi_{\rho,\beta}(r)}=\alpha,$$
where $\Psi_{\rho,\beta}(z)=z^{\beta}\lambda_{\rho,\beta}(z).$
\item[\bf b.] The radius of convexity of order $\alpha$ of $g_{\rho,\beta}$ is the smallest positive root of $$1+\frac{rg_{\rho,\beta}''(r)}{g_{\rho,\beta}'(r)}=\alpha.$$
\item[\bf c.] The radius of convexity of order $\alpha$ of $h_{\rho,\beta}$ is the smallest positive root of $$1+\frac{rh_{\rho,\beta}''(r)}{h_{\rho,\beta}'(r)}=\alpha.$$
\end{enumerate}
\end{theorem}

Now, we present some lower and upper bounds for the radii of convexity of the functions $g_{\rho,\beta}$ and $h_{\rho,\beta}$ by using the corresponding Euler-Rayleigh inequalities.

\begin{theorem}\label{th6}
For $\rho,\beta>0$ the radius of convexity $r^{c}(g_{\rho,\beta})$ of the function  $g_{\rho,\beta}$ is the smallest positive root of the equation $(zg^{\prime}_{\rho,\beta}(z))^{\prime}=0$ and satisfies the following inequalities
$$\sqrt{\frac{\Gamma(\rho+\beta)}{9\Gamma(\beta)}}<r^{c}(g_{\rho,\beta})<
\sqrt{\frac{9\Gamma(\rho+\beta)\Gamma(2\rho+\beta)}{\Delta_{81,25}(\rho,\beta)}},$$
$$\sqrt[4]{\frac{\Gamma^{2}(\rho+\beta)\Gamma(2\rho+\beta)}{\Gamma(\beta)\Delta_{81,25}(\rho,\beta)}}<r^{c}(g_{\rho,\beta})<
\sqrt{\frac{2\Gamma(\rho+\beta)\Gamma(3\rho+\beta)\Delta_{81,25}(\rho,\beta) }{\Gamma(\beta)\Gamma(3\rho+\beta)\Delta_{1458,675}(\rho,\beta)+49\Gamma^{3}(\rho+\beta)\Gamma(2\rho+\beta)}} .$$
\end{theorem}

\begin{theorem}\label{th7}
For $\rho,\beta>0$ the radius of convexity $r^{c}(h_{\rho,\beta})$ of the function $h_{\rho,\beta}$ is the smallest positive root of the equation $(zh^{\prime}_{\rho,\beta}(z))^{\prime}=0$ and satisfies the following inequalities	
$$\frac{\Gamma(\rho+\beta)}{4\Gamma(\beta)}<r^{c}(h_{\rho,\beta})<
\frac{4\Gamma(\rho+\beta)\Gamma(2\rho+\beta)}{\Delta_{16,9}(\rho,\beta)},$$
$$\sqrt{\frac{\Gamma^{2}(\rho+\beta)\Gamma(2\rho+\beta)}{\Gamma(\beta)\Delta_{16,9}(\rho,\beta)}}<r^{c}(h_{\rho,\beta})<
\frac{\Gamma(\rho+\beta)\Gamma(3\rho+\beta)\Delta_{16,9}(\rho,\beta)}{8\Gamma^{3}(\rho+\beta)\Gamma(2\rho+\beta)+
2\Gamma(\beta)\Gamma(3\rho+\beta)\Delta_{32,27}(\rho,\beta)}.$$
\end{theorem}

\subsection{Some particular cases of the main results}
It is important to mention that the Wright function is actually a generalization of a transformation of the Bessel function of the first kind. Namely, we have the relation
$$\lambda_{1,1+\nu}(z)=\phi(1,1+\nu,-z^{2})=z^{-\nu}J_{\nu}(2z)$$
where $J_{\nu}$ stands for the Bessel function of the first kind and order $\nu$. Taking into account this, it is clear that Theorem \ref{th1} in particular when $\rho=1$ and $\beta=\nu+1$ reduce to some interesting results, and one of the them naturally complements the results from \cite[Theorem 1]{Ba0}. The result on $f_{1,\nu+1}$ is new and complements \cite[Theorem 1]{Ba0}, however, the results on $g_{1,\nu+1}$ and $h_{1,\nu+1}$ are not new, they were proved in \cite[Theorem 1]{Ba0}. Thus the last two parts of Theorem \ref{th1} are natural generalizations of parts {\bf b} and {\bf c} of \cite[Theorem 1]{Ba0}.

\begin{corollary}
Let $\nu>-1$ and $\alpha\in[0,1).$
\begin{enumerate}
\item[\bf a.] The radius of starlikeness of order $\alpha$ of the function
$z \mapsto f_{1,\nu+1}(z)=\left(\Gamma(\nu+1)zJ_{\nu}(2z) \right)^{\frac{1}{\nu+1}}$
is the smallest positive root of the equation
$$2zJ^{\prime}_{\nu}(2z)+(1-\alpha(\nu+1))J_{\nu}(2z)=0.$$
\item[\bf b.] The radius of starlikeness of order $\alpha$ of the function
$z \mapsto g_{1,\nu+1}(z)=\Gamma(\nu+1)z^{1-\nu}J_{\nu}(2z)$
is the smallest positive root of the equation
$$2zJ^{\prime}_{\nu}(2z)+(1-\alpha-\nu)J_{\nu}(2z)=0.$$
\item[\bf c.] The radius of starlikeness of order $\alpha$ of the function
$z \mapsto h_{1,\nu+1}(z)=\Gamma(\nu+1)z^{1-\frac{\nu}{2}}J_{\nu}(2\sqrt{z})$
is the smallest positive root of the equation
$$2\sqrt{z}J_{\nu}'(2\sqrt{z})+(2-2\alpha-\nu)J_{\nu}(2\sqrt{z})=0.$$
\end{enumerate}
\end{corollary}

By choosing in Theorem \ref{th2} the values $\rho=1$ and $\beta=\nu+1$ we obtain the following corollary.

\begin{corollary}
If $\nu>-1$ then we have
\begin{align*}
\sqrt{\frac{\nu+1}{\nu+3}}&<r^{\star}(f_{1,\nu+1})<(\nu+1)\sqrt{\frac{(\nu+2)(\nu+3)}{\nu^3+7\nu^{2}+15\nu+13}},\\
 \sqrt[4]{\frac{(\nu+1)^{3}(\nu+2)}{\nu^3+7\nu^{2}+15\nu+13}}&<r^{\star}(f_{1,\nu+1})<
 (\nu+1)\sqrt{\frac{2(\nu+3)(\nu^3+7\nu^{2}+15\nu+13)}{\nu^5+15\nu^4+80\nu^{3}+222\nu^{2}+319\nu+196)}}.
\end{align*}
\end{corollary}

Now, by using the relation between the Wright function and the Bessel function of the first kind we can see that our main results which are given in Theorem \ref{th3} and Theorem \ref{th4} when we take $\rho=1$ and $\beta=\nu+1$ correspond to the results in \cite[Theorem 1]{aby} and \cite[Theorem 2]{aby}.

\begin{corollary}
If $\nu>-1$ then we have
\begin{align*}
\sqrt{\frac{\nu+1}{3}}&<r^{\star}(g_{1,\nu+1})<\sqrt{\frac{3(\nu+1)(\nu+2)}{4\nu+13}},\\
\sqrt[4]{\frac{(\nu+1)^{2}(\nu+2)}{4\nu+13}}&<r^{\star}(g_{1,\nu+1})<\sqrt{\frac{(\nu+1)(\nu+3)(4\nu+13)}{2(4\nu^{2}+26\nu+49)}}.
\end{align*}
\end{corollary}

 Consider the function $z \mapsto\varphi_{\nu}(z)=2^{\nu}\Gamma(\nu+1)z^{1-\nu}J_{\nu}(z),$ which is a normalized Bessel function of the first kind, considered in \cite[Theorem 1]{aby}. Since $\varphi_{\nu}(2z)=2g_{1,\nu+1}(z)$ we obtain that the above inequalities coincide with the inequalities of \cite[Theorem 1]{aby}.

\begin{corollary}
If $\nu>-1$ then we have
\begin{align*}
\frac{\nu+1}{2}&<r^{\star}(h_{1,\nu+1})<\frac{2(\nu+1)(\nu+2)}{\nu+5},\\
\frac{(\nu+1)\sqrt{\nu+2}}{\sqrt{\nu+5}}&<r^{\star}(h_{1,\nu+1})<\frac{(\nu+1)(\nu+3)(\nu+5)}{\nu^{2}+8\nu+23}.
\end{align*}
\end{corollary}

By considering that $\phi_{\nu}(4z)=4h_{1,\nu+1}(z),$ where $\Phi_{\nu}(z)=2^{\nu}\Gamma(\nu+1)z^{1-\frac{\nu}{2}}J_{\nu}(\sqrt{z}),$ we can see that the above inequalities correspond to the results of \cite[Theorem 2]{aby}.

Finally, we mention that if we take $\rho=1$, $\beta=\nu+1$ with $\nu>-1$ in Theorem \ref{th6} and Theorem \ref{th7} we can see that the following
inequalities correspond to the results which are given in \cite[Theorem 6]{abo} and \cite[Theorem 7]{abo}, respectively.

\begin{corollary}
If $\nu>-1$ then we have
\begin{align*}
\frac{\sqrt{\nu+1}}{3}&<r^{c}(g_{1,\nu+1})<3\sqrt{\frac{(\nu+1)(\nu+2)}{56\nu+137}},\\
\sqrt[4]{\frac{(\nu+1)^{2}(\nu+2)}{56\nu+137}}&<r^{c}(g_{1,\nu+1})<\sqrt{\frac{(\nu+1)(\nu+3)(56\nu+137)}{2(208\nu^{2}+1172\nu+1693)}}.
\end{align*}
\end{corollary}

\begin{corollary}
If $\nu>-1$ then we have
\begin{align*}
\frac{\nu+1}{4}&<r^{c}(h_{1,\nu+1})<\frac{4(\nu+1)(\nu+2)}{7\nu+23},\\
\sqrt{\frac{(\nu+1)^{2}(\nu+2)}{7\nu+23}}&<r^{c}(h_{1,\nu+1})<\frac{(\nu+1)(\nu+3)(7\nu+23)}{2(9\nu^{2}+60\nu+115)}.
\end{align*}
\end{corollary}

\subsection{Problems for further research} It is interesting to see how far the properties of Bessel functions of the first kind may be extended to apply to the Wright function. In this paper we can see that those properties of Bessel functions which comes from the fact that it is entire can be extended to the Wright function without a major difficulty. However, we would like to see whether other properties of the Bessel functions of the first kind can be extended or not to Wright functions. Here is a short list on possible open questions/problems, which worth to be studied:

\begin{enumerate}
\item[\bf 1.] {\em What can we say about the monotonicity of the zeros $\lambda_{\rho,\beta,n}$ with respect to $\beta$ (or $\rho$)?} The answer to this question would ensure that it would be possible to obtain necessary and sufficient conditions on the parameters $\rho$ and $\beta$ such that the normalized forms of the Wright function belong to certain class of univalent functions, like starlike, convex or spirallike. Such kind of results would improve the existing results in the literature (see \cite{mustafa,p1,p2,raza}).
\item[\bf 2.] {\em Is it possible to express for fixed $n$ the derivative of the zeros $\lambda_{\rho,\beta,n}$ with respect to $\beta$ (or $\rho$)?} In \cite{Ba1} the Watson formulae for the derivative of the zeros of Bessel function of the first kind and its derivative played an important role in obtaining necessary conditions for the order of the normalized Bessel functions of the first kind such that these functions to belong to the class of convex functions.
\item[\bf 3.] {\em Is it possible to use continued fractions to obtain the order of starlikeness and convexity of the normalized Wright functions?}
\end{enumerate}

Each of the above problems seems to be difficult to solve because the Wright function is not a solution of a second order homogeneous linear differential equation (as the Bessel function is) and although its power series structure is similar to that of Bessel functions it seems that its properties are more difficult to be studied.

\section{\bf Proofs of the main results}
\setcounter{equation}{0}

\begin{proof}[\bf Proof of Lemma \ref{lemma}]
The proof of the reality of the zeros is given in \cite{sanjeev} by using two somehow similar approaches. Now, since the growth order of the entire function $\phi(\rho,\beta,\cdot)$ is $(\rho+1)^{-1}$ (see \cite{mainardi}), which is a non-integer number and lies in $(0,1),$ it follows that indeed the Wright function has infinitely many zeros. Since the Wright function is entire its infinite product clearly exists, and in view of the Hadamard theorem on growth order of entire function it follows that its canonical representation is exactly what we have in Lemma \ref{lemma}. Using the infinite product representation we get that
\begin{equation}\label{quo1}
\frac{\Psi_{\rho,\beta}'(z)}{\Psi_{\rho,\beta}(z)}=\frac{\beta}{z}+\frac{\lambda^{\prime}_{\rho,\beta}(z)}{\lambda_{\rho,\beta}(z)}=
\frac{\beta}{z}+\sum_{n\geq1}\frac{2z}{z^{2}-\lambda^{2}_{\rho,\beta,n}}.
\end{equation}
Differentiating both sides of \eqref{quo1} we have
$$\frac{d}{dz}\left(\frac{\Psi^{\prime}_{\rho,\beta}(z)}{\Psi_{\rho,\beta}(z)}\right)=-\frac{\beta}{z^2}-2\sum_{n\geq1}\frac{z^{2}+\lambda^{2}_{\rho,\beta,n}}{(z^{2}-\lambda^{2}_{\rho,\beta,n})^{2}},\text{ \ \ } z\neq\lambda_{\rho,\beta,n}.$$
Since the expression on the right-hand side is real and negative for $z$ real and $\rho,\beta>0$, the quotient $\Psi^{\prime}_{\rho,\beta}/\Psi_{\rho,\beta}$ is a strictly decreasing function from $+\infty$ to $-\infty$ as $z$ increases through real values over the open interval $\left(\lambda_{\rho,\beta,n},\lambda_{\rho,\beta,n+1}\right),$ $n\in\mathbb{N}$. Hence, the function $\Psi^{\prime}_{\rho,\beta}$ vanishes just once between two consecutive zeros of the function $\lambda_{\rho,\beta}.$
\end{proof}

\begin{proof}[\bf Proof of Theorem \ref{th1}]
We need to show that the inequalities
\begin{equation}\label{equ3}
\real\left( \frac{zf^{\prime }(z)}{f(z)}\right) \geq\alpha,\  \real\left( \frac{zg^{\prime }(z)}{g(z)}\right) \geq\alpha \text{ \ \ and \ \ } \real\left( \frac{zh^{\prime }(z)}{h(z)}\right) \geq\alpha
\end{equation}
are valid for $z\in\mathbb{D}_{r^\star_{\alpha}}\left( f_{\rho,\beta}\right)$, $z\in\mathbb{D}_{r^\star_{\alpha}}\left( g_{\rho,\beta}\right)$ and $z\in\mathbb{D}_{r^\star_{\alpha}}\left( h_{\rho,\beta}\right)$ respectively, and each of the above inequalities does not hold in any larger disk. By definition we get
$$f_{\rho,\beta}(z)=\left(z^{\beta}\Gamma(\beta)\lambda_{\rho,\beta}(z)\right)^{\frac{1}{\beta}},$$
$$g_{\rho,\beta}(z)=z\Gamma(\beta)\lambda_{\rho,\beta}(z),$$
$$h_{\rho,\beta}(z)=z\Gamma(\beta)\lambda_{\rho,\beta}(\sqrt{z}).$$
Logarithmic derivation yields
$$\frac{zf_{\rho,\beta}^{\prime}(z)}{f_{\rho,\beta}(z)}=1+\frac{1}{\beta}\left( \frac{z\lambda^{\prime}_{\rho,\beta}(z)}{\lambda_{\rho,\beta}(z)}\right)=1-\frac{1}{\beta}\sum_{n\geq1}\frac{2z^2}{\lambda_{\rho,\beta,n}^2-z^2},$$
$$\frac{zg_{\rho,\beta}^{\prime}(z)}{g_{\rho,\beta}(z)}=1+\left( \frac{z\lambda^{\prime}_{\rho,\beta}(z)}{\lambda_{\rho,\beta}(z)}\right)=1-\sum_{n\geq1}\frac{2z^2}{\lambda_{\rho,\beta,n}^2-z^2},$$
$$\frac{zh_{\rho,\beta}^{\prime}(z)}{h_{\rho,\beta}(z)}=1+\frac{1}{2}\left(\sqrt{z} \frac{\lambda^{\prime}_{\rho,\beta}(\sqrt{z})}{\lambda_{\rho,\beta}(\sqrt{z})}\right)=1-\sum_{n\geq1}\frac{z}{\lambda_{\rho,\beta,n}^2-z}.$$
It is known \cite{Ba0} that if $z\in\mathbb{C}$ and $\theta\in\mathbb{R}$ are such that $\theta>\left| z\right| $ then
\begin{equation}\label{eqszasz}\frac{\left| z\right|}{\theta-\left| z\right|}\geq\real\left( \frac{z}{\theta-z}\right).\end{equation}
Then the inequality
$$\frac{\left| z\right|^2}{\lambda^2_{\rho,\beta,n}-\left| z\right|^2}\geq\real\left( \frac{z^2}{\lambda^2_{\rho,\beta,n}-z^2}\right)$$
is valid for every $\rho>0$, $\beta>0$, $n\in\mathbb{N}$ and $\left| z\right|<\lambda_{\rho,\beta,1}$. Therefore,
$$\real\left(\frac{zf^{\prime }_{\rho,\beta}(z)}{f_{\rho,\beta}(z)}\right)=1-\frac{1}{\beta}\real\left( \sum_{n\geq1}\frac{2z^2}{\lambda_{\rho,\beta,n}^2-z^2}\right)\geq1-\frac{1}{\beta}\sum_{n\geq1}\frac{2\left|z\right|^2 }{\lambda_{\rho,\beta,n}^2-\left|z\right|^2}=\frac{\left| z\right| f^{\prime }_{\rho,\beta}(\left| z\right| )}{f_{\rho,\beta}(\left| z\right| )},$$
$$\real\left(\frac{zg_{\rho,\beta}^{\prime}(z)}{g_{\rho,\beta}(z)}\right)=1-\real\left( \sum_{n\geq1}\frac{2z^2}{\lambda_{\rho,\beta,n}^2-z^2}\right) \geq1-\sum_{n\geq1}\frac{2\left|z\right| ^2}{\lambda_{\rho,\beta,n}^2-\left| z\right| ^2}=\frac{\left| z\right| g^{\prime }_{\rho,\beta}(\left| z\right| )}{g_{\rho,\beta}(\left| z\right| )}, $$
and
$$\real\left( \frac{zh_{\rho,\beta}^{\prime}(z)}{h_{\rho,\beta}(z)}\right)=1-\real\left( \sum_{n\geq1}\frac{z}{\lambda_{\rho,\beta,n}^2-z}\right)\geq1- \sum_{n\geq1}\frac{\left| z\right| }{\lambda_{\rho,\beta,n}^2-\left| z\right| }=\frac{\left| z\right| h^{\prime }_{\rho,\beta}(\left| z\right| )}{h_{\rho,\beta}(\left| z\right| )},  $$
where equalities are attained only when $z=\left| z\right| =r$. The latter inequalities and the minimum principle for harmonic functions imply that the corresponding inequalities in (\ref{equ3}) hold if only if $\left| z\right| < x_{\rho,\beta,1},\left| z\right| < y_{\rho,\beta,1}$ and $\left| z\right| < z_{\rho,\beta,1}$, respectively, where $x_{\rho,\beta,1}$, $y_{\rho,\beta,1}$ and $z_{\rho,\beta,1}$ are the smallest positive roots of the equations
$$
\frac{rf_{\rho,\beta}^{\prime}(r)}{f_{\rho,\beta}(r)}=\alpha,\ \frac{rg_{\rho,\beta}^{\prime}(r)}{g_{\rho,\beta}(r)}=\alpha \text{ \ \ and \ \ } \frac{rh_{\rho,\beta}^{\prime}(r)}{h_{\rho,\beta}(r)}=\alpha,
$$
which are equivalent to
$$r\lambda^{\prime}_{\rho,\beta}(r)-(\alpha-1)\beta\lambda_{\rho,\beta}(r)=0,\ \ \ r\lambda^{\prime}_{\rho,\beta}(r)-(\alpha-1)\lambda_{\rho,\beta}(r)=0$$
and $$\sqrt{r}\lambda^{\prime}_{\rho,\beta}(\sqrt{r})-(\alpha-1)\lambda_{\rho,\beta}(\sqrt{r})=0$$
In other words, we proved that
$$\inf_{z\in\mathbb{D}_{r}}\real\left( \frac{zf_{\rho,\beta}^{\prime}(z)}{f_{\rho,\beta}(z)} \right)=\frac{rf_{\rho,\beta}^{\prime}(r)}{f_{\rho,\beta}(r)}=F_{\rho,\beta}(r), \text{ \ \ \ \ } \inf_{z\in\mathbb{D}_{r}}\real\left( \frac{zg_{\rho,\beta}^{\prime}(z)}{g_{\rho,\beta}(z)} \right)=\frac{rg_{\rho,\beta}^{\prime}(r)}{g_{\rho,\beta}(r)} =G_{\rho,\beta}(r)$$
and
$$\inf_{z\in\mathbb{D}_{r}}\real\left( \frac{zh_{\rho,\beta}^{\prime}(z)}{h_{\rho,\beta}(z)} \right)=\frac{rh_{\rho,\beta}^{\prime}(r)}{h_{\rho,\beta}(r)}=H_{\rho,\beta}(r).$$
Since the real functions $F_{\rho,\beta},G_{\rho,\beta},H_{\rho,\beta}:(0,\lambda_{\rho,\beta,1})\longrightarrow\mathbb{R}$ are decreasing, and take the limits
$$\lim_{r\searrow0}F_{\rho,\beta}(r)=\lim_{r\searrow0}G_{\rho,\beta}(r)=\lim_{r\searrow0}H_{\rho,\beta}(r)=1$$
and $$\lim_{r\nearrow\lambda_{\rho,\beta,1}}F_{\rho,\beta}(r)=\lim_{r\nearrow\lambda_{\rho,\beta,1}}G_{\rho,\beta}(r)=
\lim_{r\nearrow\lambda_{\rho,\beta,1}}H_{\rho,\beta}(r)=-\infty,$$
it follows that the inequalities in (\ref{equ3}) indeed hold for $z\in\mathbb{D}_{x_{\rho,\beta,1}},$ $z\in\mathbb{D}_{y_{\rho,\beta,1}}$
and $z\in\mathbb{D}_{z_{\rho,\beta,1}},$ respectively.
\end{proof}

\begin{proof}[\bf Proof of Theorem \ref{th2}]
The radius of starlikeness of normalized Wright function$f_{\rho,\beta}$ corresponds to the radius of starlikeness of the function $\Psi_{\rho,\beta}(z)=z^{\beta}\lambda_{\rho,\beta}(z).$ The infinite series representations of the function $\Psi^{\prime}_{\rho,\beta}$ and its derivative read as follows
\begin{equation}\label{rosf1}
\Upsilon_{\rho,\beta}(z)=\Psi^{\prime}_{\rho,\beta}(z)=\sum_{n\geq0}\frac{(-1)^{n}(2n+\beta)}{n!\Gamma(n\rho+\beta)}z^{2n+\beta-1},
\end{equation}
\begin{equation}\label{rosf2}
\Upsilon^{\prime}_{\rho,\beta}(z)=\sum_{n\geq0}\frac{(-1)^{n}(2n+\beta)(2n+\beta-1)}{n!\Gamma(n\rho+\beta)}z^{2n+\beta-2}.
\end{equation}
In view of Lemma \ref{lemma} the function $z\mapsto z^{1-\beta}\Upsilon_{\rho,\beta}(z)$ belongs to the Laguerre-P\'{o}lya class $\mathcal{LP}$. Hence, the zeros of the function $\Upsilon_{\rho,\beta}$ are all real. Suppose that $\iota_{\rho,\beta,n}$'s are the positive zeros of the function $\Upsilon_{\rho,\beta}$. The expression $\Upsilon_{\rho,\beta}(z)$ can be written as
\begin{equation}\label{rosf3}
\Gamma(\beta)\Upsilon_{\rho,\beta}(z)=\beta z^{\beta-1}\prod_{n\geq1}\left( 1-\frac{z^{2}}{\iota^{2}_{\rho,\beta,n}}\right).
\end{equation}
By logarithmic derivation of both sides of \eqref{rosf3} for $|z|<\iota_{\rho,\beta,1}$we obtain
\begin{equation}\label{rosf4}
\frac{z\Upsilon^{\prime}_{\rho,\beta}(z)}{\Upsilon_{\rho,\beta}(z)}-({\beta-1})=-2\sum_{n\geq1}\frac{z^2}{\iota_{\rho,\beta,n}^2-z^2}=-2\sum_{n\geq1}\sum_{k\geq0}\frac{z^{2k+2}}{\iota_{\rho,\beta,n}^{2k+2}}=-2\sum_{k\geq0}\chi_{k+1}z^{2k+2},
\end{equation}
where $\chi_{k}=\sum_{n\geq1}\iota^{-2k}_{\rho,\beta,n}$. Thus, by using the relations \eqref{rosf1}, \eqref{rosf2} and \eqref{rosf4} we get
\begin{equation}\label{rosf5}
\frac{z\Upsilon^{\prime}_{\rho,\beta}(z)}{\Upsilon_{\rho,\beta}(z)}=\left.\sum_{n\geq0}\xi_{n}z^{2n}\right/\sum_{n\geq0}\nu_{n}z^{2n},
\end{equation}
where
$$\xi_{n}=(-1)^{n}\frac{(2n+\beta)(2n+\beta-1)}{n!\Gamma(n\rho+\beta)}\quad \text{\\ and \\}\quad\nu_{n}=(-1)^{n}\frac{(2n+\beta)}{n!\Gamma(n\rho+\beta)}.$$
By comparing the coefficients of \eqref{rosf4} and (\ref{rosf5}) we have
$$\left\{\begin{array}{l}(\beta-1)\nu_0=\xi_0\\
(\beta-1)\nu_1-2\chi_1\nu_0=\xi_1\\
(\beta-1)\nu_2-2\chi_1\nu_1-2\chi_2\nu_0=\xi_2\\
(\beta-1)\nu_3-2\chi_1\nu_2-2\chi_2\nu_1-2\chi_3\nu_0=\xi_3
\end{array}\right.,$$
which implies that
$$\chi_{1}=\frac{(\beta+2)\Gamma(\beta)}{\Gamma(\rho+\beta)},\ \chi_{2}=\frac{(\beta+2)^2}{\beta}\frac{\Gamma^2(\beta)}{\Gamma^2(\rho+\beta)}-\frac{\beta+4}{\beta}\frac{\Gamma(\beta)}{\Gamma(2\rho+\beta)}$$
and
$$\chi_{3}=\frac{(\beta+2)^3}{\beta^2}\frac{\Gamma^3(\beta)}{\Gamma^3(\rho+\beta)}-
\frac{(\beta+2)^2(\beta+4)\Gamma^2(\beta)}{2\beta^2\Gamma(\rho+\beta)\Gamma(2\rho+\beta)}+\frac{\beta+6}{2\beta}\frac{\Gamma(\beta)}{\Gamma(3\rho+\beta)}.$$
By using the Euler-Rayleigh inequalities $\chi^{-1/k}_{k}<\iota^{2}_{\rho,\beta,1}<\frac{\chi_{k}}{\chi_{k+1}}$, $k\in\{1,2\}$, we get the inequalities of the theorem.
\end{proof}

\begin{proof}[\bf Proof of Theorem \ref{th3}]
For $\alpha=0,$ in view of the second part of Theorem \ref{th1}, we have that the radius of starlikeness of order zero is the smallest positive root of the equation $(z\lambda_{\rho,\beta}(z))^{\prime}=0$. Therefore, we shall study the first positive zero of
\begin{equation}\label{equ5}
\psi_{\rho,\beta}(z)=(z\lambda_{\rho,\beta}(z))^{\prime}=\sum_{n\geq0}\frac{(-1)^{n}(2n+1)}{n!\Gamma(n\rho+\beta)}z^{2n}.
\end{equation}
We know that the function $\lambda_{\rho,\beta}$ belongs to the Laguerre-P\'{o}lya class of entire functions $\mathcal{LP}$, which is closed under differentiation. Therefore, we get that the function $\psi_{\rho,\beta}$ belongs also to the Laguerre-P\'{o}lya class. Hence, the zeros of the function $\psi_{\rho,\beta}$ are all real. Suppose that $\gamma_{\rho,\beta,n}$'s are the positive zeros of the function $\psi_{\rho,\beta}$. Then, the function $\psi_{\rho,\beta}$ has the infinite product representation as follows,
\begin{equation}\label{equ6}
\Gamma(\beta)\psi_{\rho,\beta}(z)=\prod_{n\geq1}\left( 1-\frac{z^2}{\gamma^{2}_{\rho,\beta,n}}\right),
\end{equation}
since its growth order corresponds to the growth order of the Wright function itself. If we take the logarithmic derivative of both sides of \eqref{equ6}, then for $|z|<\gamma_{\rho,\beta,1}$ we get
\begin{equation}\label{equ7}
\frac{\psi^{\prime}_{\rho,\beta}(z)}{\psi_{\rho,\beta}(z)}=\sum_{n\geq1}\frac{2z}{z^2-\gamma_{\rho,\beta,n}^2}=
-2\sum_{n\geq1}\sum_{k\geq0}\frac{z^{2k+1}}{\gamma_{\rho,\beta,n}^{2k+2}}=-2\sum_{k\geq0}\sum_{n\geq1}\frac{z^{2k+1}}{\gamma_{\rho,\beta,n}^{2k+2}}
=-2\sum_{k\geq0}\delta_{k+1}z^{2k+1},
\end{equation}
where $\delta_{k}=\sum_{n\geq1}\gamma^{-2k}_{\rho,\beta,n}.$ Moreover, in view of \eqref{equ5} we have
\begin{equation}\label{equ8}
\frac{\psi^{\prime}_{\rho,\beta}(z)}{\psi_{\rho,\beta}(z)}=-2\left.\sum_{n\geq0}a_{n}z^{2n+1}\right/\sum_{n\geq0}b_{n}z^{2n},
\end{equation}
where
$$a_{n}=\frac{(-1)^{n}(2n+3)}{n!\Gamma((n+1)\rho+\beta)} \text{ \ \and \ \ } b_{n}=\frac{(-1)^{n}(2n+1)}{n!\Gamma(n\rho+\beta)}.$$
Comparing the coefficients of \eqref{equ7} and \eqref{equ8} we obtain
$$\delta_1b_0=a_0,\ \delta_2b_0+\delta_1b_1=a_1,\ \delta_3b_0+\delta_2b_1+\delta_1b_2=a_2,$$ which yields the following Rayleigh sums
$$\delta_{1}=\frac{3\Gamma(\beta)}{\Gamma(\rho+\beta)},\ \delta_{2}=\frac{9\Gamma^2(\beta)}{\Gamma^2(\rho+\beta)}-\frac{5\Gamma(\beta)}{\Gamma(2\rho+\beta)}$$
and
$$\delta_{3}=\frac{27\Gamma^{3}(\beta)}{\Gamma^{3}(\rho+\beta)}-\frac{45}{2}\frac{\Gamma^2(\beta)}{\Gamma(\rho+\beta)\Gamma(2\rho+\beta)}+
\frac{7}{2}\frac{\Gamma(\beta)}{\Gamma(3\rho+\beta)}.$$
By using Euler-Rayleigh inequalities $\delta^{-\frac{1}{k}}_{k}<\gamma^{2}_{\rho,\beta,1}<\frac{\delta_{k}}{\delta_{k+1}}$, $k\in\{1,2\}$,
we obtain
$$\sqrt{\frac{\Gamma(\rho+\beta)}{3\Gamma(\beta)}}<r^\star(g_{\rho,\beta})<
\sqrt{\frac{3\Gamma(\rho+\beta)\Gamma(2\rho+\beta)}{\Delta_{9,5}(\rho,\beta)}},$$
$$\sqrt[4]{\frac{\Gamma^{2}(\rho+\beta)\Gamma(2\rho+\beta)}{\Gamma(\beta)\Delta_{9,5}(\rho,\beta)}}<r^\star(g_{\rho,\beta})<
\sqrt{\frac{2\Gamma(\rho+\beta)\Gamma(3\rho+\beta)\Delta_{9,5}(\rho,\beta)}{9\Gamma(\beta)\Gamma(3\rho+\beta)\Delta_{6,5}(\rho,\beta) +7\Gamma^{3}(\rho+\beta)\Gamma(2\rho+\beta)}}.$$
\end{proof}

\begin{proof}[\bf Proof of Theorem \ref{th4}]
If we take $\alpha=0$ in the third part of Theorem \ref{th1}, then we conclude that the radius of starlikeness of the function $h_{\rho,\beta}$ is actually the smallest positive root of the transcendental equation $(z\lambda_{\rho,\beta}(\sqrt{z}))^{\prime}=0$. Therefore, it is of interest to study the first positive zero of
\begin{equation}\label{rosh0}
\Omega_{\rho,\beta}(z)=(z\lambda_{\rho,\beta}(\sqrt{z}))^{\prime}=\sum_{n\geq0}\frac{(-1)^{n}(n+1)}{n!\Gamma(n\rho+\beta)}z^{n}.
\end{equation}
In view of Lemma \ref{lemma} and because of the fact that $\mathcal{LP}$ is closed under differentiation, the function $\Omega_{\rho,\beta}$ also belongs to the Laguerre-P\'{o}lya class. Assume that $\sigma_{\rho,\beta,n}$'s are the positive zeros of the function  $\Omega_{\rho,\beta}$. Thus, due to the Hadamard factorization theorem the expression $\Omega_{\rho,\beta}(z)$ can be written as
\begin{equation}\label{rosh1}
\Gamma(\beta)\Omega_{\rho,\beta}(z)=\prod_{n\geq1}\left(1-\frac{z}{\alpha_{\rho,\beta,n}}\right).
\end{equation}
By taking the logarithmic derivative of both sides of \eqref{rosh1} we have
\begin{equation}\label{rosh2}
\frac{\Omega^{\prime}_{\rho,\beta}(z)}{\Omega_{\rho,\beta}(z)}=-\sum_{k\geq0}\eta_{k+1}z^{k},\text{ \ \ \ \ } \left| z\right| <\sigma_{\rho,\beta,1},
\end{equation}
where $\eta_{k}=\sum_{n\geq1}\sigma^{-k}_{\rho,\beta,n}$. Also, by taking the derivative of \eqref{rosh0} we get
\begin{equation}\label{rosh3}
\frac{\Omega^{\prime}_{\rho,\beta}(z)}{\Omega_{\rho,\beta}(z)}=-\left.\sum_{n\geq0}c_{n}z^{n}\right/\sum_{n\geq0}d_{n}z^{n},
\end{equation}
where
$$c_{n}=\frac{(-1)^{n}(n+2)}{n!\Gamma((n+1)\rho+\beta)} \text{ \ \and \ \ }  d_{n}=\frac{(-1)^{n}(n+1)}{n!\Gamma(n\rho+\beta)}.$$
Comparing the coefficients of \eqref{rosh2} and \eqref{rosh3} we get the following Rayleigh sums
$$\eta_{1}=\frac{2\Gamma(\beta)}{\Gamma(\rho+\beta)},\ \eta_{2}=\frac{4\Gamma^2(\beta)}{\Gamma^2(\rho+\beta)}-\frac{3\Gamma(\beta)}{\Gamma(2\rho+\beta)}$$
and
$$\eta_{3}=\frac{8\Gamma^{3}(\beta)}{\Gamma^{3}(\rho+\beta)}+\frac{2\Gamma(\beta)}{\Gamma(3\rho+\beta)}-
\frac{9\Gamma^2(\beta)}{\Gamma(\rho+\beta)\Gamma(2\rho+\beta)}$$
and by using the Euler-Rayleigh inequalities $\eta^{-1/k}_{k}<\sigma_{\rho,\beta,1}<\frac{\eta_{k}}{\eta_{k+1}}$ for $k\in\{1,2\}$ we get
the inequalities of the theorem.
\end{proof}

\begin{proof}[\bf Proof of Theorem \ref{th5}]
{\bf a.} Observe that
$$1+\frac{zf^{\prime \prime}_{\rho,\beta}(z)}{f^{\prime }_{\rho,\beta}(z)}=1+\frac{z\Psi_{\rho,\beta}''(z)}{\Psi^{\prime}_{\rho,\beta}(z)}+\left( \frac{1}{\beta}-1\right) \frac{z\Psi^{\prime}_{\rho,\beta}(z)}{\Psi_{\rho,\beta}(z)}.$$
Now, we consider the following infinite product representations
$$\Gamma(\beta)\Psi_{\rho,\beta}(z)=z^{\beta}\prod_{n\geq1}\left( 1-\frac{z^{2}}{\zeta^{2}_{\rho,\beta,n}}\right),\ \Gamma(\beta)\Psi^{\prime}_{\rho,\beta}(z)=z^{\beta-1}\prod_{n\geq1}\left( 1-\frac{z^{2}}{\zeta^{\prime 2}_{\rho,\beta,n}}\right) , $$
where $\zeta_{\rho,\beta,n}$ and $\zeta^{\prime}_{\rho,\beta,n}$ are the $n$th positive roots of $\Psi_{\rho,\beta}$ and $\Psi^{\prime}_{\rho,\beta},$ respectively. Note that $\zeta_{\rho,\beta,n}$ is in fact equal to $\lambda_{\rho,\beta,n},$ however since the zeros of $\lambda_{\rho,\beta}'$ and $\Psi_{\rho,\beta}'$ do not coincide we use different notations for the zeros of the derivatives, and hence also for the zeros of $\Psi_{\rho,\beta}.$ Logarithmic differentiation on both sides of the above relations yields
$$\frac{z\Psi^{\prime}_{\rho,\beta}(z)}{\Psi_{\rho,\beta}(z)}=\beta-\sum_{n\geq1}\frac{2z^{2}}{\zeta^{2}_{\rho,\beta,n}-z^{2}},\text{ \ \ \ \ }\frac{z\Psi^{\prime \prime}_{\rho,\beta}(z)}{\Psi^{\prime}_{\rho,\beta}(z)}=\beta-1-\sum_{n\geq1}\frac{2z^{2}}{\zeta^{\prime 2}_{\rho,\beta,n}-z^{2}},$$
which implies that
$$1+\frac{zf^{\prime \prime}_{\rho,\beta}(z)}{f^{\prime }_{\rho,\beta}(z)}=1-\left( \frac{1}{\beta}-1\right) \sum_{n\geq1}\frac{2z^{2}}{\zeta^{2}_{\rho,\beta,n}-z^{2}}-\sum_{n\geq1}\frac{2z^{2}}{\zeta^{\prime 2}_{\rho,\beta,n}-z^{2}}.$$
By using the inequality \eqref{eqszasz} for $\beta\in(0,1]$ we obtain that
\begin{equation}\label{ineqconv}
\real\left(1+\frac{zf^{\prime \prime}_{\rho,\beta}(z)}{f^{\prime }_{\rho,\beta}(z)}\right)\geq 1-\left( \frac{1}{\beta}-1\right) \sum_{n\geq1}\frac{2r^{2}}{\zeta^{2}_{\rho,\beta,n}-r^{2}}-\sum_{n\geq1}\frac{2r^{2}}{\zeta^{\prime 2}_{\rho,\beta,n}-r^{2}},
\end{equation}
where $\left| z\right| =r.$ Moreover, in view of \cite[Lemma 2.1]{Ba1}, that is,
$$\alpha\real\left(\frac{z}{a-z}\right)-\real\left(\frac{z}{b-z}\right)\geq \alpha\frac{|z|}{a-|z|}-\frac{|z|}{b-|z|},$$
where $a>b>0,$ $z\in\mathbb{C}$ such that $|z|<b,$ we obtain that \eqref{ineqconv} is also valid when $\beta>1$ for all  $z\in\mathbb{D}_{\zeta^{\prime}_{\rho,\beta,1}}.$ Here we used tacitly that the zeros of $\zeta_{\rho,\beta,n}$ and $\zeta^{\prime}_{\rho,\beta,n}$ interlace according to Lemma \ref{lemma}, that is, we have $\zeta^{\prime}_{\rho,\beta,1}<\zeta_{\rho,\beta,1}.$ Now, the above deduced inequalities imply for $r\in(0,\zeta^{\prime}_{\rho,\beta,1})$
$$\inf_{z\in\mathbb{D}_{r}}\left\lbrace \real\left(1+\frac{zf^{\prime \prime}_{\rho,\beta}(z)}{f^{\prime }_{\rho,\beta}(z)}\right)\right\rbrace =1+r\frac{f^{\prime \prime}_{\rho,\beta}(r)}{f^{\prime }_{\rho,\beta}(r)}. $$
On the other hand, the function $u_{\rho,\beta}:(0,\zeta^{\prime}_{\rho,\beta,1})\rightarrow\mathbb{R},$ defined by
$$u_{\rho,\beta}(r)=1+\frac{rf^{\prime \prime}_{\rho,\beta}(r)}{f^{\prime }_{\rho,\beta}(r)},$$
is strictly decreasing when $\beta\in(0,1].$ Moreover, it is also strictly decreasing when $\beta>1$ since
\begin{align*}
u^{\prime}_{\rho,\beta}(r)&=-\left(\frac{1}{\beta}-1\right)\sum_{n\geq1}\frac{4r\zeta^{2}_{\rho,\beta,n}}{(\zeta^{2}_{\rho,\beta,n}-r^{2})^{2}}-\sum_{n\geq1}\frac{4r\zeta^{\prime 2}_{\rho,\beta,n}}{(\zeta^{\prime 2}_{\rho,\beta,n}-r^{2})^{2}}
\\&<\sum_{n\geq1}\frac{4r\zeta^{2}_{\rho,\beta,n}}{(\zeta^{2}_{\rho,\beta,n}-r^{2})^{2}}-\sum_{n\geq1}\frac{4r\zeta^{\prime 2}_{\rho,\beta,n}}{(\zeta^{\prime 2}_{\rho,\beta,n}-r^{2})^{2}}<0
\end{align*}
for $r\in(0,\zeta^{\prime}_{\rho,\beta,1}).$ Here we used the interlacing property of the zeros stated in Lemma \ref{lemma}. Observe also that $\lim_{r\searrow0}u_{\rho,\beta}(r)=1$ and $\lim_{r\nearrow\zeta^{\prime}_{\rho,\beta,1}}u_{\rho,\beta}(r)=-\infty$, which means that for $z\in\mathbb{D}_{r_{1}}$ we get
$$ \real\left( 1+\frac{zf^{\prime \prime}_{\rho,\beta}(z)}{f^{\prime }_{\rho,\beta}(z)}\right)>\alpha $$
if and only if $r_{1}$ is the unique root of
$$1+\frac{rf^{\prime \prime}_{\rho,\beta}(r)}{f^{\prime }_{\rho,\beta}(r)}=\alpha$$
situated in $(0,\zeta^{\prime}_{\rho,\beta,1}).$

{\bf b.} Since $g_{\rho,\beta}\in\mathcal{LP}$ it follows that $g_{\rho,\beta}'\in\mathcal{LP},$ and since their growth orders (which coincide according to the theory of entire functions) are equal to $(\rho+1)^{-1},$ we get via the Hadamard theorem the Weierstrassian canonical representation
$$g_{\rho,\beta}'(z)=\prod_{n\geq1}\left(1-\frac{z^2}{\vartheta_{\rho,\beta,n}^2}\right).$$ Logarithmic derivation of both sides yields
$$1+\frac{zg''_{\rho,\beta}(z)}{g'_{\rho,\beta}(z)}=1-\sum_{n\geq1}\frac{2z^2}{\vartheta_{\rho,\beta,n}^2-z^2}.$$
Application of the inequality \eqref{eqszasz} implies that
$$\real\left(1+\frac{zg''_{\rho,\beta}(z)}{g'_{\rho,\beta}(z)}\right)\geq1-\sum_{n\geq1}\frac{2r^2}{\vartheta_{\rho,\beta,n}^2-r^2},$$
where $|z|=r.$ Thus, for $r\in(0,\vartheta_{\rho,\beta,1})$ we get
$$\inf_{z\in\mathbb{D}_r}\left\{\real\left(1+\frac{zg_{\rho,\beta}''(z)}{g_{\rho,\beta}'(z)}\right)\right\}=
1-\sum_{n\geq1}\frac{2r^2}{\vartheta_{\rho,\beta,n}^2-r^2}=1+\frac{rg_{\rho,\beta}''(r)}{g_{\rho,\beta}'(r)}.$$
The function
$v_{\rho,\beta}:(0,\vartheta_{\rho,\beta,1})\to\mathbb{R},$ defined by
$$v_{\rho,\beta}(r)=1+\frac{rg_{\rho,\beta}''(r)}{g_{\rho,\beta}'(r)},$$ is strictly decreasing
and $$\lim_{r\searrow0}v_{\rho,\beta}(r)=1, \ \ \ \lim_{r\nearrow\vartheta_{\rho,\beta,1}}v_{\rho,\beta}(r)=-\infty.$$
Consequently, the equation
$$1+\frac{rg_{\rho,\beta}''(r)}{g_{\rho,\beta}'(r)}=\alpha$$  has a unique root
$r_2$ in $(0,\vartheta_{\rho,\beta,1}).$ In other words, we have
$$\real\left(1+\frac{zg_{\rho,\beta}''(z)}{g_{\rho,\beta}'(z)}\right)>\alpha, \
z\in{\mathbb{D}_{r_2}} \ \ \textrm{and} \ \
\inf_{z\in{\mathbb{D}_{r_2}}}\left\{\real\left(1+\frac{zg_{\rho,\beta}''(z)}{g_{\rho,\beta}'(z)}\right)\right\}=\alpha.$$

{\bf c.} By using again the fact that the zeros of the Wright function $\lambda_{\rho,\beta}$ are all real and in view of the Hadamard theorem we obtain
$$h_{\rho,\beta}'(z)=\prod_{n\geq1}\left(1-\frac{z}{\tau_{\rho,\beta,n}}\right),$$ which implies that
$$1+\frac{zh_{\nu}''(z)}{h_{\nu}'(z)}=1-\sum_{n\geq1}\frac{z}{\tau_{\rho,\beta,n}-z}.$$
Let $r\in\left(0,\tau_{\rho,\beta,1}\right)$ be a fixed number. The minimum principle for harmonic functions and inequality \eqref{eqszasz} imply that for $z\in\mathbb{D}_r$ we have
\begin{align*}
\real&\left(1+\frac{zh_{\rho,\beta}''(z)}{h_{\rho,\beta}'(z)}\right)=
\real\left(1-\sum_{n\geq1}\frac{z}{\tau_{\rho,\beta,n}-z}\right)\geq
\min_{|z|=r}\real\left(1-\sum_{n\geq1}\frac{z}{\tau_{\rho,\beta,n}-z}\right)\\
&=\min_{|z|=r}\left(1-\sum_{n\geq1}\real\frac{z}{\tau_{\rho,\beta,n}-z}\right)
\geq1-\sum_{n\geq1}\frac{r}{\tau_{\rho,\beta,n}-r}=1+\frac{rh_{\nu}''(r)}{h_{\nu}'(r)}.
\end{align*}
Consequently, it follows that
$$\inf_{z\in{\mathbb{D}_r}}\left\{\real\left(1+\frac{zh_{\rho,\beta}''(z)}{h_{\rho,\beta}'(z)}\right)\right\}=
1+\frac{rh_{\rho,\beta}''(r)}{h_{\rho,\beta}'(r)}.$$
Now, let $r_3$  be the smallest positive  root of the equation
\begin{equation}\label{ab}1+\frac{rh_{\rho,\beta}''(r)}{h_{\rho,\beta}'(r)}=\alpha\end{equation}
For $z\in\mathbb{D}_{r_3}$ we have
$$\real\left(1+\frac{zh_{\rho,\beta}''(z)}{h_{\rho,\beta}'(z)}\right)>\alpha.$$ In order to finish the proof, we need to show that the equation \eqref{ab} has a unique root in $\left(0,\tau_{\rho,\beta,1}\right).$
But, the equation \eqref{ab} is equivalent to
$$w_{\nu}(r)=1-\alpha-\sum_{n\geq1}\frac{r}{\tau_{\rho,\beta,n}-r}=0,$$
and we have $$\lim_{r\searrow0}w_{\nu}(r)=1-\alpha>0, \   \   \  \lim_{r\nearrow\tau_{\rho,\beta,1}}w_{\nu}(r)=-\infty.$$
Now, since the function $w_{\nu}$ is strictly decreasing on $(0,\tau_{\rho,\beta,1}),$ it follows that the equation $w_{\nu}(r)=0$ has a unique root.
\end{proof}

\begin{proof}[\bf Proof of Theorem \ref{th6}]
By using the infinite series representations of the Wright function and its derivative we obtain
$$\Theta_{\rho,\beta}(z)=(zg^{\prime}_{\rho,\beta})^{\prime}=1+\sum_{n\geq1}\frac{(-1)^{n}(2n+1)^{2}\Gamma(\beta)}{n!\Gamma(n\rho+\beta)}z^{2n}.$$
We know that the function $g_{\rho,\beta}$ belongs to the Laguerre-P\'{o}lya class and $\mathcal{LP}$ is closed under differentiation. Thus, the function $\Theta_{\rho,\beta}$ belongs also to the Laguerre-P\'{o}lya class and hence its zeros are all real. Assume that $\varsigma_{\rho,\beta,n}$'s are the positive zeros of the function $\Theta_{\rho,\beta}$. The function $\Theta_{\rho,\beta}$ can be written as follows
$$\Theta_{\rho,\beta}(z)=\prod_{n\geq1}\left( 1-\frac{z^{2}}{\varsigma^{2}_{\rho,\beta,n}}\right),$$
which for $|z|<\varsigma_{\rho,\beta,1}$ yields
\begin{equation}\label{equ13}
\frac{\Theta^{\prime}_{\rho,\beta}(z)}{\Theta_{\rho,\beta}(z)}=-2\sum_{n\geq1}\frac{z}{\varsigma^{2}_{\rho,\beta,n}-z^2}=
-2\sum_{n\geq1}\sum_{k\geq0}\frac{z^{2k+1}}{\varsigma^{2k+2}_{\rho,\beta,n}}=-2\sum_{k\geq0}\sum_{n\geq1}\frac{z^{2k+1}}{\varsigma^{2k+2}_{\rho,\beta,n}}=
-2\sum_{k\geq0}\kappa_{k+1}z^{2k+1},
\end{equation}
where $\kappa_{k}=\sum_{n\geq1}\varsigma^{-2k}_{\rho,\beta,n}$. On the other hand, we have
\begin{equation}\label{equ14}
\frac{\Theta^{\prime}_{\rho,\beta}(z)}{\Theta_{\rho,\beta}(z)}=-2\left.\sum_{n\geq0}q_{n}z^{2n+1}\right/\sum_{n\geq0}r_{n}z^{2n},
\end{equation}
where
$$q_{n}=\frac{(-1)^{n}(2n+3)^{2}\Gamma(\beta)}{n!\Gamma((n+1)\rho+\beta)}\text{ \ \and \ \ } r_{n}=\frac{(-1)^{n}(2n+1)^{2}\Gamma(\beta)}{n!\Gamma(n\rho+\beta)}.$$
By comparing the coefficients of (\ref{equ13}) and (\ref{equ14}) we obtain
$$\kappa_{1}=\frac{9\Gamma(\beta)}{\Gamma(\rho+\beta)},\ \
\kappa_{2}=\frac{81\Gamma^{2}(\beta)}{\Gamma^2(\rho+\beta)}-\frac{25\Gamma(\beta)}{\Gamma(2\rho+\beta)}$$
and
$$\kappa_{3}=\frac{729\Gamma^3(\beta)}{\Gamma^3(\rho+\beta)}+\frac{49}{2}\frac{\Gamma(\beta)}{\Gamma(3\rho+\beta)}-
\frac{675}{2}\frac{\Gamma^2(\beta)}{\Gamma(\rho+\beta)\Gamma(2\rho+\beta)}.$$
By using the Euler-Rayleigh inequalities $\kappa^{-1/k}_{k}<\varsigma^{2}_{\rho,\beta,1}<\frac{\kappa_{k}}{\kappa_{k+1}}$ for $k\in\{1,2\}$ we obtain the inequalities of the theorem.
\end{proof}

\begin{proof}[\bf Proof of Theorem \ref{th7}]
By definition we have
\begin{equation}\label{equ15}
\omega_{\rho,\beta}(z)=(zh^{\prime}_{\rho,\beta}(z))^{\prime}=1+\sum_{n\geq1}\frac{(-1)^{n}(n+1)^{2}\Gamma(\beta)}{n!\Gamma(n\rho+\beta)}^{}z^{n}.
\end{equation}
Moreover, we know that $h_{\rho,\beta}$ belongs to the Laguerre-P\'{o}lya class $\mathcal{LP},$ and consequently the function $\omega_{\rho,\beta}$ belongs also to the Laguerre-P\'{o}lya class. In other words, the zeros of the function $\omega_{\rho,\beta}$ are all real. Assume that $\varrho_{\rho,\beta,n}$'s are the positive zeros of the function $\omega_{\rho,\beta}$. In this case, the function $\omega_{\rho,\beta}$ has the infinite product representation as follows
\begin{equation}\label{equ16}
\omega_{\rho,\beta}(z)=\prod_{n\geq1}\left( 1-\frac{z}{\varrho_{\rho,\beta,n}}\right).
\end{equation}
By taking the logarithmic derivative of both sides of (\ref{equ16}) for $|z|<\varrho_{\rho,\beta,1}$ we have
\begin{equation}\label{equ17}
\frac{\omega^{\prime}_{\rho,\beta}(z)}{\omega_{\rho,\beta}(z)}=-\sum_{k\geq0}\mu_{k+1}z^{k},
\end{equation}
where $\mu_{k}=\sum_{n\geq1}\varrho^{-k}_{\rho,\beta,n}$. In addition, by using the derivative of infinite sum representation of (\ref{equ15}) we obtain
\begin{equation}\label{equ18}
\frac{\omega^{\prime}_{\rho,\beta}(z)}{\omega_{\rho,\beta}(z)}=-\left.\sum_{n\geq0}t_{n}z^{n}\right/\sum_{n\geq0}s_{n}z^{n}
\end{equation}
where
$$t_{n}=\frac{(-1)^{n}(n+2)^{2}\Gamma(\beta)}{n!\Gamma((n+1)\rho+\beta)} \text{ \ \and \ \ } s_{n}=\frac{(-1)^{n}(n+1)^{2}\Gamma(\beta)}{n!\Gamma(n\rho+\beta)}. $$
By comparing the coefficients of (\ref{equ17}) and (\ref{equ18}) we get
$$\mu_{1}=\frac{4\Gamma(\beta)}{\Gamma(\rho+\beta)},\text{ \ \ \ \ } \mu_{2}=\frac{16\Gamma^{2}(\beta)}{\Gamma^{2}(\rho+\beta)}-\frac{9\Gamma(\beta)}{\Gamma(2\rho+\beta)}$$
and
$$\mu_{3}=\frac{64\Gamma^{3}(\beta)}{\Gamma^{3}(\rho+\beta)}+\frac{8\Gamma(\beta)}{\Gamma(3\rho+\beta)}-\frac{54\Gamma^{2}(\beta)}{\Gamma(\rho+\beta)\Gamma(2\rho+\beta)}.$$
By considering the Euler-Rayleigh inequalities $\mu^{-1/k}_{k}<\varrho_{\rho,\beta,1}<\frac{\mu_{k}}{\mu_{k+1}}$,  $k\in\{1,2\},$ we have the inequalities of the theorem.
\end{proof}

\end{document}